\newcommand{\switch}[2]{#1} %for arxiv

\switch
{\documentclass[a4paper,10pt]{article}
\usepackage{paper}
\usepackage{hyperref}
\hypersetup{colorlinks=true,citecolor=black,linkcolor=black,anchorcolor=black,filecolor=black,menucolor=black,urlcolor=black
}
}%
{
\documentclass[smallextended,referee,envcountsect,]{svjour3}
\usepackage{graphicx}
\usepackage{amsmath,amssymb,pifont}
\usepackage{bm}%makes better bold fornt in math
\usepackage{manfnt}
\usepackage{enumerate}
\usepackage{wrapfig}
\usepackage[
sorting=nyt,
bibencoding=auto,
backend=biber,
maxnames=10,
autolang=other,
doi=false,
url=false,
style=numeric-comp,
isbn=false]{biblatex}
\usepackage{csquotes}
\AtEveryBibitem{\clearlist{language}}
\AtEveryBibitem{\clearlist{note}}
\renewbibmacro{in:}{}

\addbibresource{Library.bib}
\DeclareFieldFormat{postnote}{#1}
\DeclareFieldFormat{multipostnote}{#1}
\DeclareFieldFormat{pages}{#1}
\def\parit#1{\medskip\noindent{\it #1}}
\def\parbf#1{\medskip\noindent{\bf #1}}
\def\qeds{\qed\par\medskip}
\def\eps{\varepsilon}
\def\epsilon{\varepsilon}
\def\ge{\geqslant}
\def\le{\leqslant}
\def\geq{\geqslant}
\def\leq{\leqslant}
\def\df{\mathrel{{:}{=}}}
\def\dist{\operatorname{\rm dist}}
\let\oldcdot\cdot
\def\cdot{{\hskip0.5pt\oldcdot\hskip0.5pt}}
\newcommand*{\z}[1]{#1\nobreak\discretionary{}%
            {\hbox{$\mathsurround=0pt #1$}}{}}
\newcounter{thm}[section]
\renewcommand{\thethm}{\thesection.\arabic{thm}}
\def\claim#1{\par\medskip\noindent\refstepcounter{thm}\hbox{\bf #1 \thesection.\arabic{thm}.}
\it\ %\ignorespaces
}
\def\endclaim{
\par\medskip}
\newenvironment{thm}{\claim}{\endclaim}
\usepackage{hyperref}
\journalname{JOTA}
}

\def\thetitle{\switch{Cyclic projections in Hadamard spaces}{Cyclic Projections in Hadamard Spaces}}
\def\theauthors{Alexander Lytchak and Anton Petrunin}
\hypersetup{
pdftitle={\thetitle},
pdfauthor={\theauthors}
}

\switch
{

}{\institute{Alexander Lytchak \at
Karlsruhe Institute of Technology,\\
D-76128, Karlsruhe, Germany\\
alytchak@math.uni-koeln.de
\and
Anton Petrunin \at
Math. Dept. PSU, \\
University Park, PA 16802, USA.\\
petrunin@math.psu.edu
}
}

\begin{document}

\title{\thetitle}	
\author{\theauthors}
\date{}
\maketitle

\begin{abstract}
We show that cyclic products of projections onto convex subsets of Hadamard spaces can behave in a more complicated way than in Hilbert spaces, resolving a problem formulated by Miroslav Bačák.
Namely, we construct an example of convex subsets in a Hadamard space  such that the corresponding cyclic product of projections is not asymptotically regular.

\ 

\switch{}{Communicated by Aviv Gibali}
\end{abstract}

%\begin{keywords}
%CAT(0), Hadamard space, asymptotic regularity, metric projection
%\end{keywords}
%\subjclass[2010]{53C20, 53C21, 53C23}
%\pagestyle{empty}\thispagestyle{empty}

\section{Introduction}

The method of cyclic  projections is a classical algorithm seeking  an intersection point of a finite family  $C_1,\dots, C_k$ of  closed convex subsets in a Hilbert  space~$X$.
Denote by $P_i$ the closest-point projection $X\to C_i$; it sends a point $x\in X$ to the (necessarily unique) point $P_i(x)$ in $C_i$ that minimizes the distance to~$x$.
Given a point $x\in X$ consider the sequence $x_n=P^n(x)$, where
$P$ is the  cyclic composition of projections $P= P_1\circ \dots \circ P_k$.
The method of cyclic projections analyzes the sequence $(x_n)$, tries to find a limit point $x_{\infty}$, to show $x_{\infty} \in C_1\cap\z\dots\cap C_k$, and to understand the  rate of convergence.

Let us list some results in the area.
\emph{If the intersection  $C_1\cap \dots\cap C_k$ is non-empty, then $(x_n)$ always converges weakly to some point in $C_1\cap \dots\cap C_k$} \cite{bregman}.
However, this convergence does not need to be strong \cite{hundal}.
\emph{If, in addition,  $C_i$ are  linear subspaces, then the convergence is strong} \cite{vonneumann,halperin}.
If the intersection  $C_1\cap \dots\cap C_k$ is not assumed to be non-empty,  the analysis of the sequence $(x_n)$ is more complicated.
However, in \cite{bauschke} it has been established that the cyclic product $P= P_1\circ \dots \circ P_k$ is \emph{asymptotically regular};
by definition, this means, that for any starting point  $x\in X$, we have  $|x_n-x_{n+1}|\to 0$ as $n\to \infty$.
The rates of convergence, respectively, the rates of asymptotic regularity have been investigated in several works, see, for instance \cite{bauschke-borwein-lewis, Kohlenbach}.
For further reference, see \cite{bauschke-borwein-lewis, deutsch, Bac, Bac2, AFLN}.

More recently, the method of cyclic  projections has been investigated beyond the setting of Hilbert spaces in so-called Hadamard spaces (also known as CAT(0) spaces, or globally non-positively curved spaces in the sense of Alexandrov).
This class of metric spaces includes hyperbolic spaces, metric trees, as well as complete simply-connected Riemannian manifolds of non-positive curvature;
it has played an important role in many areas of mathematics  in the last decades.
We assume some familiarity with Hadamard spaces, refer the reader to \cite{ballmannnotes,ballmannbook, BBI,BH,AKP,AKP_inv} as general references on this subject.
For the introduction and applications of the method of cyclic projections in Hadamard spaces,
see \cite{bachak-book}, \cite[Section 6.8]{Bac}, and the references therein.

Hadamard spaces are defined (loosely speaking) by the property that their distance function is at least as convex as the distance function on a Hilbert space.
In particular, Hadamard spaces contain a huge variety of convex subsets;
closest point projections to closed convex subsets are well-defined and $1$-Lipschitz,
and the questions discussed above about cyclic projections are absolutely meaningful in a Hadamard space $X$.

Many results discussed above have been transferred from the linear setting of Hilbert spaces
to general Hadamard spaces.
For instance, \emph{if the subsets $C_i$ have a non-empty intersection, then 
the cyclic product of projection $P$ is asymptotically regular and, for any initial point
$x\in X$, the sequence  $x_n\z=P^n (x)$ converges \emph{weakly} to a point $x_{\infty} \in C_1\cap \dots\cap C_k$}
\cite{asymptotic, Bac2}.
(The weak topology on Hadamard spaces is discussed in \cite{Bac, bachak-book, lytchak-petrunin}).
%Moreover, if the  subsets $C_i$ have a non-empty intersection, then the cyclic product of projections $P$ is asymptotically regular;
The rate of convergence in this setting has been studied in \cite{KLN}.

Therefore it is somewhat surprising, that the fundamental result of Heinz Bauschke \cite{bauschke} for (possibly) non-intersecting convex subsets $C_i$ does not admit  a generalization to the
 setting of Hadamard spaces.  The following main result of this paper provides a negative answer to  the question of Miroslav Bačák \cite[Problem 6.13]{Bac}.

\begin{thm}{Theorem}\label{thm}
There exist a Hadamard space $X$ and compact convex subsets $C_1,\dots,C_k$ in $X$ such that the composition of the closest-point projections $P\z= P_1\circ \dots \circ P_k$ is not asymptotically regular.
\end{thm}

%If the sets $C_1,\dots,C_k$ have a common intersection, then such examples are impossible \cite{asymptotic,Bac2,Bac}.

We provide an explicit example with $X$ being a product of two trees, proving the theorem for $k=3$.
Setting $C_3=\dots=C_k$ defines examples for any $k\geq 3$.

In this example, all subsets $C_i$ are isometric to the unit interval, the projections $P_i$ map all of these segments isometrically onto $C_i$ and the composition $P\z=P_1\circ P_2\circ P_3$ maps $C_1$ to itself isometrically but exchanges the endpoints
of this interval.
A stronger version of the theorem is proved in the appendix;
it requires a somewhat deeper understanding of the geometry of Hadamard spaces.
It seems possible, but would require some non-trivial technical work, to adapt the example
from the appendix so that the Hadamard space becomes a smooth Hadamard  manifold.

On the other hand,  in the case $k=2$,  the result of Heinz Bauschke \cite{bauschke}  does admit a generalization;  in this case, the algorithm is sufficiently simple to be controlled explicitly, even providing an optimal rate of asymptotic regularity.
As it was pointed out by an anonymous referee, the following statement follows from \cite[Theorem 3.3]{asymptotic}, under the additional assumption of the existence of a fixed point of the composition $P$.

\begin{thm}{Proposition}\label{prop}
Let $C_1,C_2$ be two closed convex subsets of a Hadamard space $X$.
Then the composition $P\z=  P_1\circ P_2$ is asymptotically  regular.

Moreover, $|x_n-x_{n+1}| =o (\frac  {1} {\sqrt n})$ for any  $x\in X$ and $x_n \z= P^n (x)$.
\end{thm}

Here and further we denote by $|x-y|$ the distance between points $x$ and $y$ in any metric space, even without linear structure.

Examples given by the real axis $C_1 \subset \mathbb{R}^2$ and the set
\[C_2  = \{\,(x,y):x>0, y \geq 1+ x^{ -\epsilon}\,\}\]
reveal that the convergence rate in Proposition~\ref{prop} cannot be improved to $O (n^{-\frac 1 2  -\epsilon})$ for any $\epsilon >0$. 

This also shows that the optimal rate of asymptotic regularity for cyclic product of projections on two convex subsets is the same for the Euclidean plane and general Hadamard spaces. 

\switch{
\parbf{Acknowledgments.}
We thank Miroslav Bačák and Nina Lebedeva for helpful comments and conversations.
Let us also thank the anonymous referees for careful reading and useful suggestions.
Alexander Lytchak was partially supported by the DFG grant, no. 281071066, TRR 191.
Anton Petrunin was partially supported by the NSF grant, DMS-2005279.
}{}

\section{\switch{Three segments in a product of two tripods}{Three Segments in a Product of Two Tripods}}\label{sec:tripods}

\switch{}{In this section we prove Theorem~\ref{thm}.}

\parit{Proof\switch{ of \ref{thm}}{}.}
A union of three unit segments that share one endpoint with the induced length metric will be called a \emph{tripod}.  
Consider two tripods $S$ and $T$ and the product space $X= S\times T$.
\begin{figure}[h!]
\vskip0mm
\centering
\includegraphics{\switch{mppics/}{}Fig20}
\end{figure}
Our space $X$ is a product of two trees, thus of two Hadamard spaces.
Hence $X$ is a Hadamard space.

Denote by $a$, $b$, $c$ and $u$, $v$, $w$ the sides of $S$ and $T$ respectively.

The following diagram shows 3 isometric copies of $2{\times}2$-square in $X$; they are obtained as the products of two pairs of sides in $S$ and $T$ as labeled.

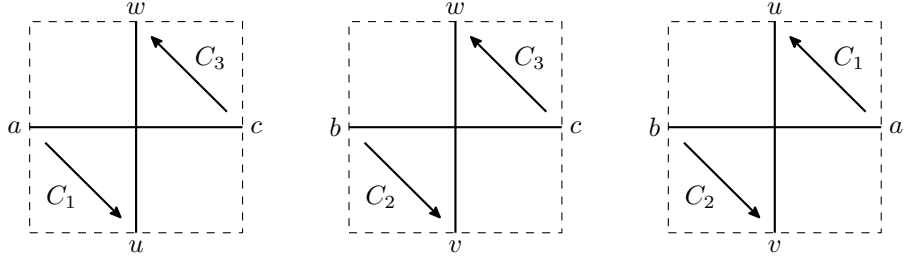
\begin{figure}[ht!]
\vskip0mm
\centering
\includegraphics{\switch{mppics/}{}Fig30}
\end{figure}

Consider the segments $C_1$, $C_2$, and $C_3$ shown on the diagram;
they all have slope $-1$ and project to each other isometrically.
Note that each projection $P_i$ reverses the shown orientation.
It follows that the composition  $P=P_1\circ P_2 \circ P_3$ sends the segment  $C_1$ to itself isometrically and changes the orientation of the segment.
In particular, $P$ exchanges the ends of the segment, hence $P$ is not asymptotically regular.
(In fact, for an end $e$ of $C_1$, and any $n$, we have $|P^n (e)\z-P^{n+1}(e)|=1$.)

Finally, setting $C_3=\dots=C_k$ defines examples for any $k\geq 3$.\qeds
  
\section{\switch{Two sets}{Two Sets}}

\switch{}{In this section we prove Proposition~\ref{prop}.}

\switch{\mbox}{}{\parit{Proof\switch{ of \ref{prop}}{}.}}
By definition, $x_n \in C_1$ for all~$n$.
Set $y_{n+1}= P_2\circ  P^n(x)$, 
so $y_1=P_2(x)$, $x_1=P_1(y_1)$, $y_2=P_2(x_1)$, and so on.
Further set 
\begin{align*}
r_n&:=|x_n-x_{n+1}|,\\
s_n&:=|y_n-y_{n+1}|.
\end{align*}

\begin{wrapfigure}[12]{r}{52mm}
\switch{\vskip-0mm}{\vskip-9mm}
\centering
\includegraphics{\switch{mppics/}{}Fig10}
\end{wrapfigure}

Since the closest-point projection is nonexpanding, we get
\[s_1 \geq r_1 \geq s_2\geq r_2\geq \dots
\switch{\eqlbl{eq:rsrsrsrs}}{\eqno(1)}
\]

Set
\begin{align*}
a_n &\df |x_n-y_n|= \dist_{C_1}y_n,\\
 b_n &\df |y_{n+1}-x_n|= \dist_{C_2}x_n.
\end{align*}
Note that
\[a_1 \geq b_1 \geq a_2 \geq b_2 \geq \dots
\switch{\eqlbl{eq:abababa}}{\eqno(2)}\]

Since $C_1$ is convex and $x_{n}\in C_1$ lies at the minimal distance from $y_n$, we have $\measuredangle[x_{n}\,{}^{x_{n-1}}_{y_n}]\ge \tfrac\pi2$. 
Since $X$ is a Hadamard space,
\[r_n^2  \leq b_{n}^2 - a_{n+1}^2.\]
Therefore, \switch{\ref{eq:abababa}}{(2)} implies that 
\[\sum_{n} r_n ^2\le b_1^2.\]
By \switch{\ref{eq:rsrsrsrs}}{(1)}, $r_n$ is non-increasing.
Therefore, $r_n = o(\tfrac1{\sqrt{n}})$.
\qeds

\switch{}{
\section{Conclusions}

We have shown that  a cyclic product of $k\ge 3$ projections onto convex subsets of a Hadamard space does not need to be asymptotically regular, even
if the convex subsets involved are compact. This should be seen in contrast to the asymptotic regularity of such maps in Hilbert spaces and to the fact that many other results
about cyclic projections generalize easily from the linear setting to the setting of Hadamard spaces.  On the other hand, we show that a cyclic product of two projections to convex subsets of Hadamard spaces  must always be  asymptotically regular.  

\begin{acknowledgements}
We thank Miroslav Bačák and Nina Lebedeva for helpful comments and conversations.
Let us also thank the anonymous referees for careful reading and useful suggestions.

Alexander Lytchak was partially supported by the DFG grant, no. 281071066, TRR 191.
Anton Petrunin was partially supported by the NSF grant, DMS-2005279.
Data sharing not applicable to this article as no datasets were generated or analysed during the current study.
%%%??? нужно для Journal of Optimization Theory and Applications
\end{acknowledgements}
}
\appendix
\switch{\section*{Appendix: Three discs} \stepcounter{section}}%
{\section*{Appendix: Three Discs} \stepcounter{section}}
\label{sec:discs}

While the cyclic product of projections $P$ constructed in Section~\ref{sec:tripods} is not asymptotically regular, its square $P^2$ is the identity on $C_1$, in particular, $P^2$ \emph{is} asymptotically regular.  The construction in Section~\ref{sec:tripods}
produces a M\"obius band $B$ divided into three rectangles and a map from $B$ to a  Hadamard space that is distance-preserving on each rectangle.

In this appendix, we produce a Hadamard space that contains an embedding of a twisted  solid torus with arbitrary twisting angle, such that the  solid torus consists of
3 isometrically embedded flat cylinders.
In this case, we obtain again 3 projections onto convex sets, each of them isometric to a  Euclidean disc, the bases of the cylinders.
Then the cyclic product of these projections is the rotation of a disc by the prescribed twisting angle $\alpha$.
In  particular, if $\tfrac \alpha \pi$ is irrational, then any power of this cyclic product of projections may not be  asymptotically regular.  

\begin{thm}{Theorem}\label{thm:powers}
There is a cyclic projection $P$ as in Theorem~\ref{thm} such that any of its power $P^m$ is not asymptotically regular.
\end{thm}

\parit{Proof of \ref{thm:powers}.}
Fix an angle $\alpha$ and a small $\epsilon>0$.
Consider the closed $\epsilon$-neighborhood $A$ of a closed geodesic $\gamma$ in the unit sphere~$\mathbb{S}^3$.
Note that the boundary of $A$ is a saddle surface in $\mathbb{S}^3$;
hence it has curvature bounded from above by 1.
Thus, $A$ is a compact Riemannian manifold with boundary, such that the curvature of the interior and of the boundary is bounded from above by 1.
Therefore, by the result of Stephanie Alexander,  David Berg, and Richard Bishop \cite {ABB-1993}, \textit{$A$ equipped with the induced intrinsic metric is locally} CAT(1).
The universal cover $\tilde A$ of $A$ with its induced metric is locally CAT(1) as well. 
Since $\tilde A$ does not contain closed geodesics, it is CAT(1), by  the generalized Hadamard--Cartan theorem
\cite[8.13.3]{AKP}, \cite[6.8+6.9]{ballmannnotes}, \cite{bowditch}.

Denote by $E$ the inverse image of $\gamma$ in $\tilde A$.
The isometry group of $\tilde A$ 
contains the group of translations along $E$ and 
the rotations that fix $E$.
Let $T$  be the composition of translation along $E$  of length $2\cdot\pi +10\cdot\epsilon$ and the rotation by angle $\alpha$.
The element $T$ generates a discrete subgroup $\Gamma$ in the group of isometries  of $\tilde A$ that acts freely and discretely on $\tilde A$.

Set $Y =\tilde A/\Gamma$.
Since $\eps$ is small, any nontrivial element of $\Gamma$ moves
every point of $\tilde A$ by more than $2\cdot\pi$.
Therefore, $Y$ is a compact locally CAT(1) space that does not contain closed geodesics of length less than $2\cdot\pi$.
Hence, by  the generalized Hadamard--Cartan theorem \cite{AKP}, $Y$ is CAT(1).
By construction, $Y$ is locally isometric to $\mathbb{S}^3$ outside its boundary $B$.
The projection of $E$ to $Y$ is a closed geodesic $G$ of length $2\cdot\pi +10\cdot\epsilon$.

Denote by $X$ the Euclidean cone over $Y$;
since $Y$ is CAT(1), we get that $X$ is CAT(0); see \cite{AKP}.
Moreover, $X$ is locally Euclidean outside its \emph{boundary} --- the cone over $B$.

The cone $Z$ over the closed geodesic $G$ is the Euclidean cone over a circle of length $2\cdot\pi +10\cdot\epsilon$.
By construction, $Z$ is a locally convex subset of $X$.
Hence, $Z$ is a convex subset of~$X$ \cite[2.2.12]{AKP_inv}.
Let us consider a geodesic triangle
$[q_1q_2q_3]$
in $Z$ that surrounds the origin $o$ of the cone $Z$.

By construction, the sides of the triangle $[q_1q_2q_3]$ lie in the flat part of $X$.
Thus, we can find a small $r>0$ such that the $2\cdot r$-neighborhood $U_1$ of the geodesic $[q_1q_2]$ is isometric to a convex subset of the Euclidean space.
We can assume that $2\cdot r$-neighborhoods $U_2$ of $[q_2q_3]$ and $U_3$ of $[q_3q_1]$ have the same property.

Denote by $C_i$ the disc of radius $r$ centered at $q_i$ and being orthogonal to~$Z$.
By construction, $C_i$ and $C_{i+1}$, for $i=1,2,3\pmod 3$   are contained in~$U_i$.
Since $Z$ is convex, $C_i$ and $C_{i+1}$ are \emph{parallel} inside $U_i$, thus their convex hull  $Q_i$ is isometric to the  cylinder $C_i \times [q_i,q_{i+1}]$ with bottom and top $C_i$ and $C_{i+1}$.
In particular, the projection $P_i$ defines an isometry $C_{i+1}\to C_{i}$.

By construction, the composition $P=P_1\circ P_2\circ P_3\colon C_1\to C_1$ rotates $C_1$ by angle $\alpha$.
If $\tfrac\alpha\pi$ is irrational, then $P$, as well as all its powers, are \emph{not} asymptotically regular.

As before, setting $C_3=\dots=C_k$ defines examples for any $k\geq 3$.
\qeds

{\sloppy
\printbibliography[heading=bibintoc]
\fussy
}

\switch{%\Addresses
}{}
\end{document}